\documentclass[11pt]{amsart}
\theoremstyle{plain}
\newtheorem{thm}{Theorem}[section]
\newtheorem{theorem}[thm]{Theorem}

\newtheorem{lemma}[thm]{Lemma}
\newtheorem{corollary}[thm]{Corollary}
\newtheorem{proposition}[thm]{Proposition}
\theoremstyle{definition}
\newtheorem{remark}[thm]{Remark}

\newtheorem{definition}[thm]{Definition}

\numberwithin{equation}{section}

\newcommand{\sH}{{\mathcal H}}

\newcommand{\sO}{{\mathcal O}}

\newcommand{\C}{{\mathbb C}}

\newcommand{\BP}{{\mathbb P}}
\newcommand{\Q}{{\mathbb Q}}
\newcommand{\R}{{\mathbb R}}

\newcommand{\Z}{{\mathbb Z}}

\title {Local structure of principally polarized stable Lagrangian fibrations   }

\author{Jun-Muk Hwang, Keiji Oguiso}

\address{Jun-Muk Hwang, Korea Institute for Advanced Study, Hoegiro 87, Seoul, 130-722, Korea} \email{jmhwang@kias.re.kr}

\address{Keiji Oguiso, Department of Mathematics, Osaka University\\
Toyonaka 560-0043 Osaka, Japan and  Korea Institute for Advanced Study, Hoegiro 87, Seoul, 130-722, Korea} \email{oguiso@math.sci.osaka-u.ac.jp}

\thanks{Jun-Muk Hwang is supported
by National Researcher Program 2010-0020413 of NRF and MEST, and Keiji Oguiso
is supported by JSPS Program 22340009 and by KIAS Scholar Program}

\begin{document}

\maketitle

\begin{abstract}
A holomorphic Lagrangian fibration is stable if the characteristic
cycles of the singular fibers are of type $I_m, 1 \leq m <\infty,$
or $A_{\infty}$. We will give a complete description of the local
structure of a stable Lagrangian fibration when it is principally
polarized. In particular, we give an explicit form of the period
map of such a fibration and conversely, for a  period map of the
described type, we construct a principally polarized stable
Lagrangian fibration with the given period map. This enables us to
give a number of examples exhibiting interesting behavior of the
characteristic cycles.
\end{abstract}

\section{Introduction}
For a holomorphic symplectic manifold $(M, \omega)$, i.e., a
$2n$-dimensional complex manifold with a holomorphic symplectic
form $\omega \in H^0(M, \Omega^2_M)$, a proper flat morphism $f:M
\to B$ over an $n$-dimensional complex manifold $B$ is called a
(holomorphic) Lagrangian fibration if all smooth fibers are
Lagrangian submanifolds of $M$. The discriminant $D \subset B$,
i.e., the set of critical  values of $f$, is a hypersurface if it
is non-empty. In
\cite{HO1}, the structure of the singular fiber of $f$ at a
general point $b\in D$ was studied. By introducing the notion of
characteristic cycles, \cite{HO1} shows that the structure of
such a singular fiber can be described in a manner completely
parallel to Kodaira's classification (\cite{Kd}, see also V. 7 in \cite{BHPV})
 of singular fibers of elliptic fibrations. Furthermore, to study the multiplicity of the
 singular fibers, \cite{HO2} generalized the stable reduction theory of elliptic fibrations (cf. V.10 in \cite{BHPV}), explicitly describing how
 arbitrary singular fiber over a general point of $D$ can  be transformed to a stable singular fiber, a singular
fiber of particularly simple type. These results exhibit that the theory of general singular fibers of
a holomorphic Lagrangian fibration gives a very natural generalization of Kodaira's theory of elliptic fibrations.

The current work is yet another manifestation of this principle.
An important part of Kodaira's theory is the study of the
asymptotic behavior of the elliptic modular function of a given
elliptic fibration near a singular fiber. As a generalization of
this we will study the asymptotic behavior of the periods of the
abelian fibers near a general singular fiber of a holomorphic
Lagrangian fibration. Here we need to make two additional
assumptions on the Lagrangian fibration.

First, we will assume that the  singular fibers are of stable
type, i.e., its characteristic cycles are of type $I_k, 0 \leq k
\leq \infty$ ($I_{\infty}$ meaning $A_{\infty}$). As explained above, any general singular fiber can
be transformed into this form by the stable reduction (\cite{HO2}, Section 4).

The second assumption we will make is that the Lagrangian
fibration is {\em principally polarized}, in the sense explained
in Definition \ref{d.pp}. This condition is satisfied if there
exists an $f$-ample line bundle on $M \setminus f^{-1}(D)$
whose restriction on smooth
fibers give principal polarizations on the abelian varieties. This
assumption is rather restrictive compared with the setting of
\cite{HO1} and \cite{HO2}, where the only assumption was that the
fibers of $f$ are of Fujiki class.

We believe
that understanding the structure of Lagrangian fibration under these assumptions is essential for the study of
general cases. Since this special case already requires
substantial care and already provides many interesting examples
(see Section 5), we restrict our
attention to it in this paper and leave the general cases to
future study.

The main result of this paper is the following.

\begin{theorem}\label{t.1}
Let $f:M \to B$ be a principally polarized stable Lagrangian
fibration  (cf. Definition \ref{d.stable} and Definition
\ref{d.pp}). Then at a general point $b \in D$ of the
discriminant, there exists a coordinate system $(z_1, \ldots,
z_n)$ with $D$ defined by $z_n =0$, such that for a suitable
choice of an integral frame of the local system $R^1 f_* \Z$ on
$B\setminus D$, the period matrices have the form
$$\theta^i_j= \frac{\partial^2 \Psi}{\partial z_i \partial z_j}
\mbox{ for } (i,j) \neq (n,n) \mbox{ and } \theta^n_n=
\frac{\partial^2 \Psi}{\partial z_n \partial z_n} + \frac{\ell}{2
\pi \sqrt{-1}} \log z_n$$ where $\Psi$ is a holomorphic function
in $z_1, \ldots, z_n$. Here $\ell$ is the number of irreducible
components of a general singular fiber.

Conversely, given any germ of holomorphic function $\Psi(z_1,
\ldots, z_n)$ such that ${\rm Im}(\theta^i_j) >0$, there exists a
principally polarized stable Lagrangian fibration whose period
matrices are of the above form.
\end{theorem}

That the period matrix of Lagrangian fibration is the Hessian of a
potential function is a well-known consequence of the action-angle
variables (cf. \cite{DM}). The logarithmic behavior of the
multi-valued part reflects the stability assumption on the
singular fiber. The novelty in Theorem \ref{t.1} lies in the
choice of the variable $z_n$ through which these two aspects are
intertwined. The existence of $z_n$  follows from  the fact proved
in Proposition \ref{p.coordi} that the characteristic foliation
accounts for the degenerate
 part of the polarization restricted to the fixed part of the monodromy.
 The  proof of this uses a version of the Monodromy Theorem from the theory of the degeneration of Hodge structures
 and the topological property of the stable singular fiber.

The converse direction in Theorem \ref{t.1} is shown by explicitly
constructing a principally polarized stable Lagrangian fibration
from a given potential function $\Psi(z)$ imaginary part of whose
Hessian matrix is positive definite. This part is a sort of
generalization of Nakamura's construction \cite{Na} of toroidal
degeneration of principally polarized abelian varieties over
$1$-dimensional small disk. Using our construction, we shall give
a concrete $4$-dimensional example of principally polarized stable
Lagrangian fibration in which the types of characteristic cycles
of singular fibers change fiber by fiber, too. To our knowledge,
such an example has not been noticed previously. In fact, most of
the previous constructions of singular fibers of Lagrangian
fibrations have used product construction from elliptic
fibrations.

\section{Stable Lagrangian fibrations}

\begin{definition}\label{d.stable} A {\em Lagrangian fibration} is a proper flat morphism $f:M
\to B$ from a holomorphic symplectic manifold $(M, \omega)$ of
dimension $2n$ to a complex manifold $B$
 of dimension $n$ such that the smooth locus of each fiber is a Lagrangian submanifold of $M$.
 The {\em discriminant} $D \subset B$ is
  the set of the critical values of $f$, which is a hypersurface in $B$ if it is non-empty. {\em Throughout this paper, we assume that $D$ is non-empty.}
We say that $f$ is a
 {\em stable Lagrangian fibration} if $D \subset B$ is a submanifold and
 each singular fiber $f^{-1}(b), b\in D,$ is stable, i.e., it is
reduced and the
characteristic cycle in the sense of \cite{HO1} is of type $I_{k},
1 \leq k \leq \infty.$ By the description in \cite{HO1}, this is
equivalent to saying that  $f^{-1}(b)$ is reduced and its normalization is a
disjoint union of a finite number of compact complex manifolds
$Y^1, \ldots, Y^{\ell}$ for some positive integer $\ell$ such that \begin{enumerate} \item[(i)]
each $Y^i$ is a $\BP^1$-bundle over an $(n-1)$-dimensional complex torus $A^i$ whose
fibers  are sent to characteristic leaves of $f^{-1}(b)$ in the sense of
\cite{HO1}, i.e., for a defining function $h \in \sO(B)$ of the
divisor $D$, the Hamiltonian vector field $\iota_{\omega}(f^* dh)$,
where $\iota_{\omega}: \Omega^1_M \to T(M)$ is the vector bundle
isomorphism induced by $\omega$, is tangent to the image of the
fibers in $M$; \item[(ii)] there exist submanifolds $S^i_1, S^i_2
\subset Y^i,$ with $S^i_1 \neq S^i_2$ except possibly when
$\ell=1,2$, such that $S^i_1 \cup S^i_2$ is a $2$-to-1 unramified
cover of $A^i$ under the $\BP^1$-bundle projection; \item[(iii)]
the normalization $\nu: \bigcup Y^i \to f^{-1}(b)$ is obtained by
the identification via a collection of biholomorphic morphisms
$g_{i}: S^i_2 \to S^{i+1}_1$ for $1 \leq i \leq \ell-1$ and
$g_{\ell}: S^{\ell}_2 \to S^1_1$ with the additional requirement $g_1 =
g_2^{-1}$ if $S^1_1 = S^1_2$ and $S^2_1 = S^2_2$  for $\ell =2$.
\end{enumerate} \end{definition}

A maximal connected union of the $\BP^1$-fibers in (i) under the
identification in (iii) is called a {\em characteristic cycle}. A
characteristic cycle can be either of finite type ($I_m$-type, $1 \leq m
<\infty$) or of infinite type $A_{\infty}$, which we also denote by
$I_{\infty}$.

 Recall (cf. \cite{HO2}
Section 4) that in a neighborhood of a general singular fiber, any
Lagrangian fibration whose fibers are of Fujiki class can be
transformed to a stable Lagrangian fibration by certain explicitly
given bimeromorphic modifications and branched covering. We will
be interested in the local property of the fibration at a point of
$D$. Thus we will make the following

({\em Assumption})   $D \subset B$ is the germ of a smooth hypersurface in an $n$-dimensional complex manifold and the fundamental group $\pi_1(B\setminus D)$ is cyclic.

\medskip
The following is immediate from  Proposition 2.2 of \cite{HO1}.

\begin{proposition}\label{p.action} Given a stable Lagrangian fibration,  we can assume that there exists an action of the complex
Lie group $\C^{n-1}$ on $M$ preserving the fibers and the
symplectic form such that $S^i_1, S^i_2$ are orbits of this action
for all $1 \leq i \leq \ell$.  This action of $\C^{n-1}$ on $Y^i$
descends to the translation action on $A^i$. The patching
biholomorphisms $g_i$ in Definition \ref{d.stable} (iii) as well
as the $\BP^1$-bundle structure in (i) are equivariant under this
action. In particular, if $S^1_1 = S^1_2$ (resp. $S^2_1=S^2_2$),
the Galois action of the double cover $S^1_1\to A^1$ (resp. $S^2_1
\to A^2$) is by a translation on the torus $S^1_1$ (resp.
$S^2_1$).
\end{proposition}

Regarding the topology of the singular fiber $f^{-1}(b), b \in D$,
we have the following.

\begin{proposition}\label{p.topology} In the setting of Definition
\ref{d.stable}, fix a component $Y^1$ of the normalization of
$f^{-1}(b)$ and set
$$Y_o:= Y^1 \setminus (S^1_1 \cup S^1_2),$$ which is equipped with a
$\C^*$-bundle structure $\varrho: Y_o \to A$ over a complex torus
$A$ of dimension $n-1$ coming from  Definition
\ref{d.stable} (i). Then there exists a (not necessarily holomorphic) continuous map $\mu: f^{-1}(b) \to
A'$ to a complex torus $A'$ isogenous to $A$ such that when $j:
Y_0 \to f^{-1}(b)$ is the natural inclusion and $\rho: Y_0 \to A
\to A'$ is the composition of $\varrho$ and an isogeny, $\mu \circ
j$ is homotopic to $\rho$.
\end{proposition}

\begin{proof}
Let us use the notation introduced in Definition \ref{d.stable}
(iii) for the description of the normalization morphism $\nu:
\bigcup Y^i \to f^{-1}(b)$.

First, we consider the case $S^i_1 \neq S^i_2$ for all $1 \leq i
\leq \ell$. Define $\widetilde{f^{-1}(b)}$ as the variety obtained
from $\bigcup Y^i$ with all the patching identification $g_{1},
\ldots, g_{\ell -1}$ such that the normalization factors through
$$\nu: \bigcup Y^i \to \widetilde{f^{-1}(b)} \to f^{-1}(b)$$ with
the second arrow given by the identification via $g_{\ell}.$ When
$\ell =1$, $\widetilde{f^{-1}(b)} = Y^1.$ The $\C^{n-1}$-action of
Proposition \ref{p.action} lifts to a $\C^{n-1}$-action on
$\widetilde{f^{-1}(b)}$. The connected unions of the images of the
$\BP^1$-fibers of $Y^i$ define finite chains of
quasi-transversally intersecting $\BP^1$'s in
$\widetilde{f^{-1}(b)},$ which we call characteristic chains. Each
characteristic chain intersects each $S^i_1$ (resp. $S^i_2$), $1
\leq i \leq \ell$ at exactly one point, inducing a morphism $
\widetilde{f^{-1}(b)}  \to A^{''}$ to a complex torus of dimension
$n-1$ biholomorphic to $A^i$'s.  This determines a biholomorphism
$\zeta: S^{\ell}_2 \to S^1_1$. Fix a point $\alpha \in
S^{\ell}_{2}$ and let $\beta = \zeta(\alpha) \in S^{1}_{1}.$ For
$t\in [0,1] \subset \R$, let $\tau_t: S^{1}_{1} \to S^{1}_{1}$ be
the translation by $t(\beta-g_{\ell}(\alpha)).$

 Define a new family of biholomorphic morphisms
$g^t_{\ell}: S^{\ell}_{2} \to S^{1}_{1}$ by $g^t_{\ell} = \tau_t
\circ g_{\ell}.$ Clearly, $g_{\ell}^0 = g_{\ell}.$ We claim that $g_{\ell}^1 = \zeta$. In fact, $\zeta^{-1} \circ g_{\ell}^1$ is an automorphism of $S^{\ell}_2$ which fixes the
point $\alpha.$ But both $g_{\ell}^1$ and $\zeta$ must be equivariant under the $\C^{n-1}$-action of
Proposition \ref{p.action}. Thus $\zeta^{-1} \circ g_{\ell}^1$ must be the identity map of $S^{\ell}_2,$
proving the claim.

 Let $f^{-1}(b)^t$ be the variety obtained from
$\widetilde{f^{-1}(b)}$ by identifying $S^{1}_{1}$ and
$S^{\ell}_{2}$ via $g^t_{\ell}$. Then $$f^{-1}(b)^0 = f^{-1}(b)$$
and $f^{-1}(b)^1$ is homeomorphic to $f^{-1}(b)$. The
$\C^{n-1}$-action descends to $f^{-1}(b)^t$ for each $t$ as
$\tau_t$ commutes with the $\C^{n-1}$-action. By abuse of
terminology, we call the maximal connected unions of the images in
$f^{-1}(b)^t$ of the characteristic chains as characteristic
cycles of $f^{-1}(b)^t$. By the $\C^{n-1}$-action, we know that
all characteristic cycles in $f^{-1}(b)^t$ are isomorphic. By our
choice of $\alpha$ and $\beta$, there exists one finite
characteristic cycle in $f^{-1}(b)^1$. Thus we get a morphism
$\mu': f^{-1}(b)^1 \to A'$ to some complex torus $A'$ isogenous to
$A^{''}$. Define $\mu: f^{-1}(b) \to A'$ as the composition of
$\mu'$ with the homeomorphism $f^{-1}(b) \to f^{-1}(b)^1$. It certainly satisfies the required property.

\medskip
Now consider the case when $\ell =1$ and $S^1_1= S^1_2$. Set $\widetilde{f^{-1}(b)} = Y^1$ and define $\zeta: S^1_2 \to S^1_1$ as the Galois
action of the double covering $S^1_1 \to A^1$ in Definition \ref{d.stable} (ii).
Then the same argument as in the previous case applies.

\medskip
Finally, consider the case when $\ell=2$, $S^1_1=S^1_2$ and
$S^2_1=S^2_2$. By Proposition \ref{p.action}, the Galois action on
$S^1_1$ (resp. $S^2_1$ ) of the double cover over $A^1$ (resp.
$A^2$) is given by a translation, say, by $\gamma_1 \in \C^{n-1}$
(resp. $\gamma_2 \in \C^{n-1}$). By the  equivariance of $g_1=
g_2^{-1}$, for each $\alpha \in S^1_1$, we have $g_1(\gamma_1
\cdot \alpha) = \gamma_2 \cdot g_1(\alpha)$. Thus by the
normalization morphism $\nu: Y^1 \cup Y^2 \to f^{-1}(b)$,  a point
$\alpha \in S^1_1$ is identified with $g_1(\alpha) \in S^2_1$, and
the point $\gamma_1 \cdot \alpha \in S^1_1$, which lies in the
$\BP^1$-fiber through $\alpha$, is identified with $\gamma_2 \cdot
g_1(\alpha)$, which lies in the $\BP^1$-fiber through
$g_1(\alpha)$. Thus we get a morphism $\mu: f^{-1}(b) \to A'$ to
an $(n-1)$-dimensional torus $A'$ whose fiber is a union of two
$\BP^1$'s identified at two points. This $\mu$ satisfies the
required property.
\end{proof}

We have the generalization of the classical action-angle correspondence as follows.

\begin{proposition}\label{p.aa} Given a stable Lagrangian
fibration $f:M \to B$, choose a Lagrangian section $\Sigma \subset
M$ of $f$.
  Then we have a natural surjective unramified  morphism $ \Phi: T^*B \to M \setminus E$
  where $E$ is the union of the irreducible components of the fibers of $f$  disjoint from $\Sigma$
  such that \begin{itemize} \item[(1)] $f \circ \Phi$ agrees with the natural projection $g: T^*B \to B$,
  \item[(2)] $\Phi$ sends the zero section of $T^*B$ to $\Sigma$ and \item[(3)] $\Phi^* \omega$ coincides
  with the standard symplectic form on $T^*B$. \end{itemize} In particular,  $\Gamma:= \Phi^{-1}(\Sigma)$
  is a Lagrangian submanifold (with many connected components) in
  $T^*B$. For each $b \in B,$  $\Phi_b := \Phi|_{T^*_b(B)}: T^*_b(B) \to f^{-1}(b)\setminus E$ is the universal covering
  and $\Gamma_b:= \Gamma \cap T^*_b(B)$ is naturally isomorphic to $H_1(f^{-1}(b) \setminus E, \Z)$.
   \end{proposition}

\begin{proof}
Over $B\setminus D$, this is just a holomorphic version (cf.
Proposition 3.5 in \cite{Hw}) of the classical action-angle
correspondence as described in Section 44 of \cite{GS}. The
statement over $D$ follows by the same argument as for the smooth
fibers. In fact, for each $b\in B$, the vector group $T^*_b(B)$
acts on the fiber $f^{-1}(b)$ with $n$-dimensional orbits on the
smooth locus of $f^{-1}(b)$ (cf. Proposition 3.3 in \cite{Hw}).
The morphism $\Phi_b$ is defined by taking the orbit map of the
point $\Sigma \cap f^{-1}(b)$ under this action, which is a
universal covering map for the smooth locus of the component of
$f^{-1}(b)$ containing $\Sigma \cap f^{-1}(b)$. This shows (1) and
(2).  The proof of (3) is the same as that of Theorem 44.2 of
\cite{GS}.
\end{proof}

\begin{proposition}\label{p.fixed} Let $b \in D$.
In the notation of Proposition \ref{p.aa}, we can assume that
    the connected component of $\Gamma$ containing each point of  $ \Gamma \cap T^*_b(B)$ is a Lagrangian section of $T^*(B) \to B$, i.e., a closed 1-form on $B$.
Let $\Gamma' \subset \Gamma$ be the union of such sections of
$\Gamma$ over $B$.  Then for each $s \in B \setminus D$,
$\Gamma'_s := \Gamma' \cap T^*_s(B)$   is a sublattice of
$\Gamma_s$ satisfying $\Gamma_s/\Gamma'_s \cong \Z$.
\end{proposition}

\begin{proof}
 Since $ f^{-1}(b) \setminus E$ is a $\C^*$-bundle over an $(n-1)$-dimensional torus,
 we see that $\Gamma \cap T^*_b(B)$ has rank $2n-1$. Thus
 $\Gamma'_s$ has rank $2n-1$. It remains to show that
 $\Gamma_s/\Gamma'_s$ is torsion-free. Suppose it has a $k$-torsion, $0< k \in \Z$, i.e.,
 there exists a point $\alpha \in \Gamma_s \setminus \Gamma'_s$ such that $k \alpha \in \Gamma'_s.$
Let $\widetilde{k\alpha}$ be a closed 1-form given by the
component of $\Gamma'$ containing $k \alpha$. Then the closed
1-form $\tilde{\alpha} = \frac{1}{k} \widetilde{k\alpha}$ is also
a component of $\Gamma'$ containing $\alpha$,  which implies
$\alpha \in \Gamma'_s$, a contradiction.
 \end{proof}

\begin{proposition}\label{p.Upsilon} In the notation of Proposition \ref{p.aa},
let $Y_o$ be the fiber of $M\setminus E$ at a point $b\in D$.  Let
$\Phi_b: T^*_b(B) \to Y_o$ be the universal covering map and
$\varrho: Y_o \to A$ be the $\C^*$-bundle over an
$(n-1)$-dimensional torus. Let $\Upsilon \subset \Gamma_b=
\Gamma'_b$ be the rank-1 sublattice corresponding to the kernel of
$$\varrho_*: H_1(Y_0, \Z) \to H_1(A, \Z).$$ Then for any $v \in
\Gamma'_b \setminus \Upsilon$, there exists $\varpi \in
H^1(f^{-1}(b), \Z)$ such that $\langle \varpi,j_*v \rangle \neq 0$ where $j_*:
H_1(Y_o, \Z) \to H_1(f^{-1}(b), \Z)$ is induced by the inclusion
$j: Y_o \subset f^{-1}(b)$.
\end{proposition}

\begin{proof}
Since $\varrho_*(v) \in H_1(A, \Z)$ is non-zero, there exists
$\varphi \in H^1(A, \Z)$ such that $\langle \varphi, \varrho_*(v) \rangle \neq
0$. Let $\varpi= \mu^* \varphi$ where the map $\mu: f^{-1}(b) \to
A$ is as defined in Proposition \ref{p.topology} satisfying $\mu
\circ j = \varrho.$ Then
$$\langle \varpi, j_*(v) \rangle  = \langle \mu^* \varphi, j_*(v) \rangle
= \langle j^* \mu^*
\varphi, v \rangle = \langle \varrho^*\varphi, v \rangle = \langle \varphi, \varrho_*(v) \rangle \neq 0.$$ \end{proof}

For a stable Lagrangian fibration $f:M \to B$,
 denote by $\Lambda$ the local system on $B\setminus
D$ defined by the lattice $\Lambda_s := H_1(f^{-1}(s), \Z)$ for $s
\in B \setminus D.$

\begin{proposition}\label{p.fixed2}
For a stable Lagrangian fibration $f:M\to B$ and $s \in B \setminus D$, fix a generator
of the cyclic
fundamental group of $\pi_1( B \setminus D,s)$ and denote by  $\tau_s: \Lambda_s \to \Lambda_s$
 the monodromy operator of the generator. Then the fixed part $\Lambda'_s \subset
\Lambda_s$ of $\tau_s$ at $s
\in B \setminus D$ has corank 1.
\end{proposition}

\begin{proof}
For any $s \in B \setminus D$, we can identify each fiber
$\Lambda_s = H_1(M_s, \Z)$ with the fiber $\Gamma_s$ of
Proposition \ref{p.aa}. Thus the result follows from Proposition
\ref{p.fixed}. \end{proof}

\section{Principally polarized stable Lagrangian fibration}

\begin{definition}\label{d.pp} Let $\Lambda$ be as in Proposition
\ref{p.fixed2}.  A {\em principal polarization} on a stable
Lagrangian fibration $f:M \to B$ is a unimodular anti-symmetric
form $Q: \wedge^2 \Lambda \to \Z_{B\setminus D}$ where
$\Z_{B\setminus D}$ denotes the constant sheaf of integers on $B
\setminus D$, which induces a principal polarization on each
smooth fibers of $f$. A stable Lagrangian fibration with a choice
of principal polarization is called a {\em principally polarized
stable Lagrangian fibration}.
\end{definition}

\begin{remark} In Definition \ref{d.pp}, the polarization on
$M \setminus f^{-1}(D)$ may not extend to an $f$-ample class
of the whole
$M$. In fact,
$f$ need not be projective. This definition is useful because
there are many situations where the polarization exists {\em a
priori} only on the smooth fibers, e.g., in Kodaira's study of
elliptic fibrations and also in our construction in Section 5.
\end{remark}

\begin{proposition}\label{p.double} Let $f: M \to B$ be a principally polarized stable Lagrangian
fibration.  Then the monodromy operator in Proposition \ref{p.fixed2} satisfies $\tau_s
\neq {\rm Id}$ and $\tau_s \circ \tau_s \neq {\rm Id}.$
\end{proposition}

\begin{proof}
If $\tau_s = {\rm Id}$, then we see that $f$ is a smooth
fibration, as  in the proof of Proposition 3.2 in \cite{Hw}. In
fact, since there is no monodromy and $f$ is polarized over
$B \setminus D$, we can extend the period map of
the abelian family on $B\setminus D$ to the whole
$B$ (\cite{Gr}, Theorem 9.5). Thus, we obtain a smooth abelian fibration
$f':M' \to B$ such that $f$ and $f'$ are bimeromorphic outside
$D$. Since $M'$ contains no rational curves and both $M$ and $M'$
have trivial canonical bundles, this implies $M$ and $M'$ are
biholomorphic, a contradiction to the non-emptiness of the
discriminant $D$ of $f$.

If $\tau_s \circ \tau_s = {\rm Id},$ take a double cover $g: B'
\to B$ branched along $D$ and let $D'= g^{-1}(D)$. Denote by
$\hat{f}: \hat{M} \to B'$ the fiber product of $f$ and $g$, which
has no monodromy on $B' \setminus D'.$ By the $\C^{n-1}$-action of
Proposition \ref{p.action} which lifts to $\hat{M}$, the following
property of $\hat{M}$ can be seen from the corresponding
properties in the case of $n=1$ (cf. Proof of Proposition 9.2 in
\cite{BHPV}): $\hat{M}$ is normal, Gorenstein  with singularities
of type $A_1 \times (\mbox{ germ of $(2n-1)$-dimensional manifold
})$ and has trivial canonical bundle. Thus we have a crepant
resolution $f': M' \to B'$, which is a family with trivial
canonical bundle and no monodromy. Then we get a contradiction as
in the previous case.
\end{proof}

\begin{lemma}\label{l.monodromy}
Let $\tau: \Lambda \to \Lambda$ be an automorphism of a lattice
such that $\Lambda':= \{ v \in \Lambda, \tau(v)=v\}$ is a
sublattice of corank 1, i.e., $\Lambda/\Lambda' \cong \Z$. If
$\tau \circ \tau \neq {\rm Id}$, then $\eta:= \tau - {\rm Id}$
satisfies $\eta \circ \eta =0.$
\end{lemma}

\begin{proof}
Note that $\Lambda'\subset {\rm Ker}(\eta)$.  The induced
automrophism $\bar{\tau}: \Lambda/\Lambda' \to \Lambda/\Lambda'$
is either ${\rm Id}$ or $-{\rm Id}$.
If $\bar{\tau} = {\rm Id}$, then for a non-zero $v \in \Lambda
\setminus \Lambda'$, we have $\tau(v) = v + \lambda$ for some
$\lambda \in \Lambda'$. Then $\eta(v) = \lambda \in \Lambda'
\subset {\rm Ker}(\eta)$. This proves that $\eta \circ \eta = 0$.
If $\bar{\tau} = - {\rm Id}$, then for a non-zero $v \in \Lambda
\setminus \Lambda'$, we have $\tau(v) = -v + \lambda$ for some
$\lambda \in \Lambda'$. Then $$\tau\circ \tau(v) = -\tau(v) +
\tau(\lambda) = -(-v + \lambda) + \lambda = v.$$ Thus $\tau\circ
\tau = {\rm Id},$ a contradiction. \end{proof}

\begin{proposition}\label{p.monodromy}
In the setting of Proposition \ref{p.double}, let $\eta:= \tau_s - {\rm Id}$. Then for any $\beta \in
{\rm Im}(\eta)$ and an element $\varphi \in H^1(M,\Z)$,
$$ \langle i^*\varphi, \beta \rangle  = \langle \varphi, i_*\beta \rangle = 0$$ where $i^*:
H^1(M, \Z) \to H^1(M_s,\Z)$ and $i_*: H_1(M_s, \Z) \to H_1(M, \Z)$
are the homomorphisms induced by the inclusion $i: M_s:=f^{-1}(s)
\subset M$.
\end{proposition}

\begin{proof}
Let $\sH$ be the local system  on $B\setminus D$  given by
$H^1(M_s, \Z), t \in B\setminus D$. Denote by $\tau^*: \sH_s
\to \sH_s$ the transformation dual to $\tau$, i.e.,
 for any
$\varpi\in H^1(M_s, \Z)$ and $u \in H_1(M_s, \Z)$,
$$ \langle \tau^*(\varpi), u \rangle  = \langle \varpi, \tau(u)\rangle.$$
By Proposition \ref{p.fixed2}, Proposition \ref{p.double} and Lemma \ref{l.monodromy},
we have $\eta \neq 0$ and $\eta\circ \eta =0$, i.e.,
$$0 \neq {\rm Im}(\eta) \subset {\rm Ker}(\eta) = \Lambda'_s.$$
 Similarly, $\eta^* := \tau^* - {\rm Id}$ is
an endomorphism of $\sH_s$ with $\eta^* \neq 0$ and $\eta^* \circ \eta^* = 0$. Since
 $i^*\varphi \in {\rm
Ker}(\eta^*)$  by (the easy half of) the global invariant cycles
theorem (cf. Theorem 4.24 of \cite{Vo}), for any $\psi \in {\rm Ker}(\eta^*)$ and $u \in \Lambda_s$,
$$ \langle \psi, \eta(u) \rangle  = \langle \eta^*(\psi), u \rangle = 0.$$ It follows that
$\langle i^*\varphi, {\rm Im}(\eta) \rangle = 0.$
\end{proof}

\begin{remark} If the family $f:M \to B$ is
projective, we could have used the Monodromy Theorem (cf. Theorem
3.15 in \cite{Vo}) in place of Proposition 3.3 and Lemma
\ref{l.monodromy} in the above proof. We have used the above
approach  because we do not want to assume that $f$ is projective.
\end{remark}

\begin{proposition}\label{p.Xi}
 For a principally polarized stable Lagrangian fibration $f:M
\to B$ and $s \in B \setminus D$, let $\tau_s: \Lambda_s \to \Lambda_s$ be the monodromy
operator of Proposition \ref{p.fixed2}, which should preserve the
 polarization $Q_s: \wedge^2 \Lambda_s \to \Z$. Setting $\eta:= \tau_s - {\rm Id}$ as in Proposition
\ref{p.monodromy}, we have $\Lambda'_s= {\rm Ker}(\eta)$.    Then ${\rm
Im}(\eta) \subset \Lambda_s$  is contained in
$$\Xi_s := \{ v\in \Lambda'_s\, \vert\, Q(v,w) =0 \mbox{ for all } w \in
\Lambda'_s\}. $$ \end{proposition}

\begin{proof} Since
$\tau_s$ preserves the polarization $Q_s$ and $\eta \circ \eta = 0$ by Lemma \ref{l.monodromy},
$$Q_s(\eta(v), u) + Q_s(v, \eta(u)) = 0 \mbox{ for all } v, u \in
\Lambda_s.$$ Thus for any $v \in \Lambda_s$ and  $u \in {\rm
Ker}(\eta)= \Lambda'_s$, we have $Q_s(\eta(v),u) = - Q_s(v, \eta(u)) = 0,$ which means
$\eta(v) \in \Xi_s.$ \end{proof}

\begin{definition}\label{d.symplectic}
Let $\Lambda$ be a free abelian group of  rank $2n$. Given a
unimodular non-degenerate anti-symmetric form $Q: \wedge^2 \Lambda
\to \Z$, a basis $\{ p_1, \ldots, p_n, q_1, \ldots, q_n\}$ of
$\Lambda$ is called a {\em symplectic basis} of $\Lambda$ with
respect to $Q$ if, in terms of the dual basis $\{ p^1, \ldots,
p^n, q^1, \ldots, q^n \}$ of ${\rm Hom}(\Lambda, \Z)$,
$$Q =  p^1 \wedge q^1 + p^2 \wedge q^2 + \cdots +  p^n \wedge
q^n.$$
\end{definition}

\begin{lemma}\label{l.1} In the setting of Definition \ref{d.symplectic}, let $\tau: \Lambda \to \Lambda$ be a group automorphism preserving $Q$.
  Assume that the subgroup $\Lambda' \subset \Lambda$ of elements fixed under $\tau$  has corank 1.
 Then there exists a symplectic basis $\{ p_1, \ldots, p_n, q_1,
\ldots, q_n\}$ such that $\{p_1, \ldots, p_n, q_1, \ldots,
q_{n-1}\} \subset \Lambda'.$
\end{lemma}

\begin{proof}
Fix a symplectic basis $\{a_1, \ldots, a_n, b_1, \ldots, b_n\}$
such that
$$ Q = a^1 \wedge b^1 + \cdots + a^n \wedge b^n.$$
The anti-symmetric form $Q|_{\Lambda'}$ must have a kernel of rank
1, i.e., $$\Xi := \{ v \in \Lambda', \; Q(v, u) = 0 \mbox{ for all
} u \in \Lambda'\}$$ has rank 1. Pick a generator $p_1$ of $\Xi$. Since $\Xi$ is primitive, i.e.,
$\Lambda/\Xi$ has no torsion,
we can write $$p_1= \alpha_1 a_1 + \cdots + \alpha_n a_n + \beta_1
b_1 + \cdots + \beta_n b_n$$ with some integers $\alpha_i,
\beta_i$ satisfying $ gcd(\alpha_1, \ldots, \alpha_n, \beta_1,
\ldots, \beta_n) = 1.$ Thus there exists integers $\alpha'_1,
\ldots, \alpha'_n, \beta'_1, \ldots, \beta'_n$ such that
$$\alpha'_1 \cdot \alpha_1 + \cdots + \alpha'_n \cdot \alpha_n +
\beta'_1 \cdot \beta_1 + \cdots + \beta'_n \cdot \beta_n = 1.$$
Let $$q_1:= -\beta'_1 a_1 - \cdots - \beta'_n a_n + \alpha'_1 b_1
+ \cdots + \alpha'_n b_n.$$ Then $Q(p_1, q_1) = 1.$ Define
$$\Lambda^{''}:= \{ v \in \Lambda, Q(p_1, v) = 0 = Q(q_1, v) \}.$$
Then $\Lambda^{''} \subset \Lambda'$ is a lattice of rank $2n-2$
such that $Q|_{\Lambda^{''}}$ is unimodular and non-degenerate (cf. \cite{GH}, the proof of Lemma in p.304). Let $\{ p_2, \ldots, p_n, q_2, \ldots, q_n\}$ be a
symplectic basis of $\Lambda^{''}$. Then $\{p_1, \ldots, p_n, q_1,
\ldots, q_n\}$ is a symplectic basis of $\Lambda$ with the
required property.
\end{proof}

\begin{proposition}\label{p.basis}
In the setting of Proposition \ref{p.fixed}, identify $\Lambda_s= H_1(M_s, \Z)$ with $\Gamma_s$ for
 $s \in B \setminus D$ as in the proof of Proposition \ref{p.fixed2}. Assume that we have a
principal polarization $Q$. Then we can find a collection of
components $\{p_1, \ldots, p_n, q_1, \ldots, q_{n-1}\}$ of
 $\Gamma'$ such that for each $x \in B \setminus D$, there exists
 $q_{n,x} \in \Gamma_x$ such that $\{ p_{1,x}, \ldots, p_{n,x}, q_{1,x}, \ldots,
 q_{n,x}\}$ is a symplectic basis of $\Lambda_x= \Gamma_x$ with respect to
 $Q_x$. \end{proposition}

 \begin{proof}
Fix a point $s \in B \setminus D$. The monodromy $\tau_s:
\Lambda_s \to \Lambda_s$ preserves the polarization $Q_s$ on $\Lambda_s$ and fixes
$\Lambda'_s= \Gamma'_s$.
Applying Lemma \ref{l.1}, we have a symplectic basis $\{ p_{1,s},
\ldots, p_{n,s}, q_{1,s}, \ldots, q_{n,s}\}$ with $p_{1,s},
\ldots, p_{n,s}, q_{1,s}, \ldots, q_{n-1,s} \in \Lambda'_s.$
Since $\Gamma'$ consists of
sections of $g : T^*B \to B$, the vectors $p_{1,s}, \ldots, p_{n,s}, q_{1,s},
\ldots, q_{n-1,s}$ uniquely determine components $p_1, \ldots,
p_n, q_1, \ldots, q_{n-1}$ of $\Lambda'$. To check the existence
of $q_{n,x}$ for any $x \in B \setminus D,$ just pick $q_{n,x}$ as
any vector in $\Lambda_x$ contained in the component of $\Lambda$
containing $q_{n,s}.$
\end{proof}

\begin{proposition}\label{p.Upsilon2}
In the notation of Proposition \ref{p.basis}, when $b \in D$, the
vector $p_{n,b}\in \Gamma'_b$ regarded as an element of $H_1(Y_o, \Z)$ in the notation of
Proposition \ref{p.Upsilon},
 lies in the lattice $\Upsilon$ of Proposition
\ref{p.Upsilon}. \end{proposition}

\begin{proof} Suppose not. By our (Assumption) after Definition \ref{d.stable}, we may assume that
$M$ is topologically retractable to $f^{-1}(b)$ and identify $H^{1}(f^{-1}(b), \Z)$ with $ H^{1}(M, \Z)$.
 Then by Proposition \ref{p.Upsilon},
there exists $\varpi \in H^{1}(f^{-1}(b), \Z) = H^{1}(M, \Z)$ such
that $\langle \varpi, j_*p_{n,b} \rangle \neq 0$. For a point $s \in B\setminus
D$,  the choice in Proposition \ref{p.basis} implies that $p_{n,s}
\in \Xi_s$ of Proposition \ref{p.Xi}. Denote by $\varpi_s$  the
element in $H^1(M_s, \Z)$ induced by $\varpi \in H^1(M, \Z)$
 under the identification $H^{1}(f^{-1}(b), \Z) = H^{1}(M, \Z).$
Since  $j_*p_{n,b} \in H_{1}(f^{-1}(b), \Z) = H_{1}(M, \Z)$ and
the image of $p_{n,s} \in H_1(M_s, \Z)$ in $H_1(M, \Z) $ belongs to
the same class, Proposition \ref{p.monodromy} and Proposition
\ref{p.Xi} say that
$$ \langle \varpi, j_*p_{n,b} \rangle = \langle \varpi_s, p_{n,s} \rangle = 0.$$ This is a
contradiction.
\end{proof}

\begin{proposition}\label{p.dh} In Proposition \ref{p.Upsilon2}, the $\C$-linear span of $\Upsilon$ in $T^*_b(B)$ is
exactly $\C \cdot dh$ where $h \in \sO(B)$ is a defining equation
of the divisor $D$. \end{proposition}

\begin{proof}
From the definition of $\Upsilon$ in Proposition \ref{p.Upsilon},
the linear span of $\Upsilon$ is sent to a fiber of the
$\C^*$-bundle. By Definition \ref{d.stable} (i), this fiber is a
leaf of the characteristic foliation, which is given by the
Hamiltonian vector field $\iota_{\omega}(f^* dh)$ on $M$. Under
the symplecto-morphism $\Phi$ in Proposition \ref{p.aa}, this
corresponds to $\C \cdot dh$.
\end{proof}

\begin{proposition}\label{p.coordi}  Let $\{p_1, \ldots, p_n, q_1, \ldots,
q_{n-1}\}$ be as in Proposition \ref{p.basis}. Then there exists a
holomorphic coordinate system $\{ z_1, \ldots, z_n\}$ on $B$ such
that, regarded as sections of $T^*(B)$,   $$p_1 = d z_1, \;
\ldots, \; p_n = d z_n$$ and $D$ is given by $z_n =0.$
\end{proposition}

\begin{proof} Since $p_1, \ldots, p_n$ are closed 1-forms which are point-wise linearly
independent at every point of $B$, we can find coordinates $ z_1,
\ldots, z_n$ with $p_i = dz_i$. By Proposition \ref{p.dh}, we may
choose $z_n$ to be a defining equation of $D$. \end{proof}

Let us recall the  classical Riemann condition (e.g.  \cite{GH}, p.306).

\begin{proposition}\label{p.classic}
Let $V$ be a complex vector space of dimension $n$ and let
$\Lambda \subset V$ be a lattice of rank $2n$ such that
$V/\Lambda$ is an abelian variety with a principal polarization.
 For a symplectic basis $\{ p_1, \ldots, p_n, q_1,
\ldots, q_n \}$ of $\Lambda$ with respect to the principal polarization $Q: \wedge^2 \Lambda \to \Z$,  $\{p_1, \ldots,
p_n\}$ becomes a $\C$-basis of $V$ and the period matrix $
(\theta_i^j)$ defined by $$q_i = \sum_{j=1}^n \theta_i^j p_j \;\; \mbox{ in $V$} $$ is
symmetric in $(i,j)$ and ${\rm
Im}(\theta_i^j) >0.$
\end{proposition}

\begin{theorem}\label{t.main} Given a principally polarized stable
Lagrangian fibration $f:M \to B$ with a Lagrangian section $\Sigma
\subset M$,  there exists a holomorphic coordinate system $(z_1,
\ldots, z_n)$ on $B$ such that \begin{itemize} \item[(i)] $z_n=0$
is a local defining equation of $D$; \item[(ii)] on $B\setminus
D$, $dz_1, \ldots, dz_{n-1}, dz_n$ belong to $\Gamma'$ in the
notation of Proposition \ref{p.fixed}; \item[(iii)] there exists a
symplectic basis $\{p_{1,s}, \ldots, p_{n,s}, q_{1,s}, \ldots,
q_{1,n}\}$ on each $\Lambda_s= \Gamma_s, s\in B\setminus D$
satisfying $$p_{1,s}= (dz_1)_s, \ldots, p_{n,s} = (dz_n)_s$$ and
the associated period matrix
 in the sense of Proposition \ref{p.classic} is given by $$ \theta^j_i = \frac{\partial^2
\Psi}{\partial z_i
\partial z_j} + \frac{\ell}{2 \pi \sqrt{-1}} \log z_n$$
for some holomorphic function $\Psi$ on $B$, which we call a {\em
potential function} of the Lagrangian fibration,  and some integer
$\ell$.
\end{itemize}
\end{theorem}

\begin{proof}
Let $\{ p_1, \ldots, p_n, q_1, \ldots, q_{n-1}\}$ be as in
Proposition \ref{p.basis} and
Proposition \ref{p.coordi}. At a point $s \in B\setminus D$, we add
$q_{n,s}$ to get a symplectic basis of $\Lambda_s$. By analytic
continuation, we get a multi-valued 1-form $q_n$ over $B\setminus D$
such that any choice of a value $q_{n,t}$ of $q_n$ at a point $t
\in B \setminus D$, together with  $p_{1,t}, \ldots, p_{n,t},
q_{1,t}, \ldots, q_{n-1,t},$ gives a symplectic basis of
$\Lambda_t$. Using the coordinate system in Proposition
\ref{p.coordi}, we can write
$$q_i = \sum_{j=1}^n \theta^j_i dz_j, $$ where $\theta^j_i$ is a
(univalent) holomorphic function on $B$ for each $1 \leq i \leq
n-1$ and $1\leq j \leq n,$  while $\theta^j_n$ is a multi-valued
holomorphic function on $B \setminus D$ for each $1 \leq j \leq
n$. By Proposition \ref{p.classic}, $\theta^i_j= \theta^j_i$ for
each $1\leq i,j\leq n$. It follows that $\theta^j_n$ is univalent
holomorphic function on $B$ for each $1 \leq j \leq n-1$. By the choice of
$p_{n,s} \in \Xi_s$ and Proposition \ref{p.Xi},  the monodromy operator $\tau_s: \Lambda_s \to
\Lambda_s$ is of the form
$$\tau_s(q_{n,s}) = q_{n,s} + \ell p_{n,s}$$
for some integer $\ell$. Thus $$\tilde{\theta}^n_n \; := \;  \theta^n_n - \frac{\ell}{2\pi
\sqrt{-1}} \log z_n$$ is also univalent. Set $\tilde{\theta}_i^j := \theta_i^j$ if $(i,j) \neq (n,n)$.   Then
$\tilde{\theta}_i^j$ is a univalent holomorphic function on $B$ for all values of $1 \leq i,j \leq n$ and
$$ q_i = \sum_{j=1}^n \tilde{\theta}^j_i dz_j \mbox{ for } 1 \leq i \leq
n-1$$ $$ q_n = \sum_{j=1}^n \tilde{\theta}^j_n d z_j + \frac{\ell}{2 \pi
\sqrt{-1}} \log z_n \; d z_n.$$ Since $q_i$'s are closed 1-forms on
$B$, we have $$\frac{\partial \tilde{\theta}^j_i}{\partial z_k} =
\frac{\partial \tilde{\theta}^k_i}{\partial z_j} = \frac{\partial
\tilde{\theta}^i_j}{\partial z_k}$$ for any $1 \leq i,j,k \leq n$. By
Poincar\'e's lemma, there exists a holomorphic function $\Psi$
such that
$$\tilde{\theta}^j_i = \frac{\partial^2 \Psi}{\partial z_i
\partial z_j}.$$ \end{proof}

\section{Construction of principally polarized stable Lagrangian fibrations
with given potential functions}
\par
\noindent

In this section, for a sufficiently small $n$-dimensional polydisk
$B$ with coordinate $(z_1, \ldots, z_n)$, we shall construct a
principally polarized stable Lagrangian fibration $f : (M,
\omega_M) \to B$ with a given potential function $\Psi(z)$. Our
construction closely follows Nakamura's toroidal construction
[Na]. However, main differences are the following:

(i) the base space $B$ is of dimension $n$ (rather than $1$).

(ii) the total space should be not only smooth but also symplectic.

\medskip
{\it (I) Construction of a non-proper Lagrangian fibration $\tilde{M} \to B$.}

For each integer $k \in \Z$, let $E_k$ be a copy of $\C \times \C$
equipped with linear coordinates $(x_k, y_k)$. We define a complex
manifold $E$ by identifying points in $\cup_{k \in \Z} E_k$ by the
following rule: a point $(x_k, y_k)$ of $E_k$ with $x_k \neq 0$
and $y_k \neq 0$ is identified with a point $(x_{k+1}, y_{k+1})$
of $E_{k+1}$ with $x_{k+1} \neq 0$ and $y_{k+1} \neq 0$, if and
only if
$$ x_{k+1} = x_k^2 y_k, \; \mbox{ and } \; y_{k+1} =
\frac{1}{x_k}.$$ On $E$, $z_n := x_k y_k$ is a well-defined
holomorphic function independent of $k$ and $$w_n := x_k^{-k+1}
y_k^{-k}$$ is a meromorphic function independent of $k$, with
zeros and poles supported on $$\cup_{k \in \Z}(x_k y_k =0).$$
Moreover, the 2-forms $dy_k \wedge dx_k$ glue together yielding a
holomorphic symplectic form $\omega_E$ on $E$, satisfying
$$\omega_E = dz_n \wedge \frac{d w_n}{w_n}.$$

Fix coordinates
$$(z_1, \ldots, z_{n-1}, w_1, \ldots, w_{n-1})$$ on  $\C^{n-1} \times
\C^{n-1}$ and regard them as functions on  the open subset
$\C^{n-1} \times (\C^{\times})^{n-1}$ defined by $$w_1 \neq 0, \;
\ldots, \; w_{n-1} \neq 0.$$ Define $$\tilde{X} := \C^{n-1} \times
(\C^{\times})^{n-1} \times E.$$ On $\tilde{X}$, we have  the
holomorphic functions $z_1, \ldots, z_n, w_1, \ldots, w_{n-1}$ and
the meromorphic function $w_n$.

Define a morphism $\tilde{p}: \tilde{X} \to \C^n$ by $(z^1,
\ldots, z^n)$. The fiber of $\tilde{p}$ over $b$ with $z_n(b)
\not= 0$ is isomorphic to
$$(\C^{\times})^{n-1} \times \C^{\times}$$ with coordinates
$(w_1, \ldots , w_{n-1}, w_n)$ and the fiber over $b$ with $z_n(b) = 0$
is isomorphic to
$$(\C^{\times})^{n-1} \times \cup_{k \in \Z} \BP^1_k$$
where $\BP^1_k$ is a copy of the projective line $\BP^1$ with affine coordinate $y_k$.

We have a holomorphic symplectic $2$-form
$$\omega_{\tilde{X}} :=  \sum_{i=1}^{n-1} dz_i \wedge \frac{dw_i}{w_i} + \omega_E  = \sum_{i=1}^{n} dz_i \wedge \frac{dw_i}{w_i}$$
on $\tilde{X}$. From now, we regard $\tilde{X}$ as a symplectic
manifold by this symplectic form. From the coordinate expression
of  $\omega_{\tilde{X}}$ and $\tilde{p}$, it is immediate that
$\tilde{p}$ is a {\it non-proper} Lagrangian fibration.

We denote
$$\tilde{M} = \tilde{X} \times_{\C^n} B$$
where
$$B = \{(z_1, \ldots , z_{n-1}, z_n) \, \vert \, \vert z_i \vert < \epsilon (\forall i)\}$$
and $\epsilon$ is a sufficiently small positive real number. We denote the natural projection $\tilde{M} \to B$
induced from $\tilde{p}$ by
$$\tilde{f} : \tilde{M} \to B\,\, .$$
Note that the restriction $\omega_{\tilde{M}}$ of
$\omega_{\tilde{X}}$ is a symplectic $2$-form on $\tilde{M}$ and
$\tilde{f}$ is a non-proper Lagrangian fibration.

\medskip
{\it (II) Group action of $\Gamma = \Z^n$ on $\tilde{M}$.}
\medskip

Let $\Psi(z_1, z_2, \cdots , z_n)$ be a holomorphic function on $B$
such that the imaginary part ${\rm Im}\ \tilde{\theta} (z)$
of the Hessian matrix
$$\tilde{\theta} (z) =
(\frac{\partial^2 \Psi}{\partial z_i\partial z_j})$$
is positive definite and $\ell$ be a positive integer. We define
the period matrix $\theta(z)$ by
$$\theta(z) = \tilde{\theta}(z) +
\frac{\log z_n}{2\pi\sqrt{-1}}\left(\begin{array}{rr}
O_{n-1}& 0\\
0 & \ell\end{array} \right)\,\, .$$
We will write
$$\tilde{\theta} (z) = \left(\begin{array}{rr}
\tilde{\Theta}_1 (z)& \tilde{\Theta}_2(z)\\
 \tilde{\Theta}_2^t(z) & \tilde{\theta}_n^n(z)\end{array} \right)\,\, ,$$
where $\tilde{\Theta}_1 (z)$ is the $(n-1) \times (n-1)$ matrix,
$\tilde{\Theta}_2(z)$ is the $(n-1) \times 1$ matrix, $\tilde{\Theta}_2^t(z)$
is the transpose of $\tilde{\Theta}_2(z)$ and $\tilde{\theta}_n^n(z)$
is $1 \times 1$ matrix.

Set $\Gamma = \Z^{n-1} \oplus \Z$. We define a group action of $\Gamma$ on $\tilde{M}$ as follows.
 Let $\gamma = (j, m) \in \Gamma$. Then the action $T_{\gamma}: \tilde{M} \to \tilde{M}$ is defined in
 terms of the coordinate functions on $\C^{n-1} \times (\C^{\times})^{n-1} \times E_k \subset \tilde{M}$ by
$$T_{\gamma}^* z_i = z_i\,\, {\rm for}\,\, i=1, \ldots, n-1$$
$$T_{\gamma}^*(\Pi_{i=1}^{n-1} w_i^{b_i}) = {\rm exp}\,(2\pi \sqrt{-1}(j \tilde{\Theta}_1(z) b + m\tilde{\Theta}_2^t(z)b) \Pi_{i=1}^{n-1} w_i^{b_i}$$
where $b =(b_i)_{i=1}^{n-1}$ is $(n-1) \times 1$ matrix, and
$$T_{\gamma}^*x_k =  ({\rm exp}\,(2\pi \sqrt{-1}(j \tilde{\Theta}_2(z) + m\tilde{\theta}_n^n(z)))^{-1}x_{k-m\ell}$$
$$T_{\gamma}^*y_k =  {\rm exp}\,(2\pi \sqrt{-1}(j \tilde{\Theta}_2(z) + m\tilde{\theta}_n^n(z))y_{k-m\ell}\,\,
.$$ It is immediate that  $T_{\gamma}^*T_{\gamma'}^* = T_{\gamma +
\gamma'}^*$. Then $T_{\gamma} \in {\rm Aut}\,(\tilde{X}/\C^n)$ and
$\gamma \mapsto T_{\gamma}$ defines an injective group
homomorphism from $\Gamma$ to ${\rm Aut}\,(\tilde{M}/\C^n)$. Here
${\rm Aut}\,(\tilde{M}/\C^n)$ is the group of automorphisms of
$\tilde{M}$ over $\C^n$, i.e., the group of automorphisms $g$ of
$\tilde{M}$ such that $\tilde{f} \circ g = \tilde{f}$.

\begin{proposition}\label{discontinuous} The action $\Gamma$ on $\tilde{M}$
is properly discontinuous, free and symplectic, in the sense that
$T_{\gamma}^* \omega_{\tilde{M}} = \omega_{\tilde{M}}$ for each
$\gamma \in \Gamma$.
\end{proposition}

\begin{proof} Freeness of the action is clear from the description of
the action. The proof of proper discontinuity is essentially the
same as the proof of \cite{Na}, Theorem 2.6. This can be also seen
from the concrete description of fibers below in (III), at least
fiberwisely.

Let us show that the action is symplectic, i.e.,
$\omega_{\tilde{M}} = T_{\gamma}^*
\omega_{\tilde{M}}$
for each $\gamma = (j, m)$. This is a new part not considered by \cite{Na}. We have $T_{\gamma}^*dz_i = dz_i$,
$T_{\gamma}^*dw_i = {\rm exp}\, (2\pi \sqrt{-1} f_i(z)) w_i$ for all $i$, where, in terms of  the standard basis $\langle e_i \rangle_{i=1}^{n-1}$
  of $\C^{n-1},$
 $$f_i(z) = j \tilde{\Theta}_1(z)e_i + m\tilde{\Theta}_2(z)e_i$$
for $1 \le i \le n-1$ and
$$f_n(z) = j \tilde{\Theta}_2(z) + m\tilde{\theta}_n^n(z)\,\, .$$
Thus, for $i$ with $1 \le i \le n$, we have
$$T_{\gamma}^*dz_i = dz_i\,\, ,$$
$$T_{\gamma}^*\frac{dw_i}{w_i} = T_{\gamma}^*(d\log w_i)$$
$$= d\log(T_{\gamma}^*w_i) = d(2\pi\sqrt{-1}f_i(z)) + d(\log w_i)$$
$$= \frac{dw_i}{w_i} + 2\pi\sqrt{-1} \sum_{k=1}^{n} \frac{\partial f_i}{\partial z_k} dz_k\,\, .$$
Using these identities, we can compute
$$T_{\gamma}^* \omega_{\tilde{M}} = \sum_{k=1}^n T_{\gamma}^*(dz_i) \wedge
T_{\gamma}^*(\frac{dw_i}{w_i})$$
$$= \omega_{\tilde{M}} - 2\pi\sqrt{-1} \sum_{i=1}^{n}\sum_{k=1}^{n} \frac{\partial f_i}{\partial z_k} dz_k \wedge dz_i$$
$$= \omega_{\tilde{M}} - 2\pi\sqrt{-1} \sum_{1 \le i < k \le n} (\frac{\partial f_i}{\partial z_k} - \frac{\partial f_k}{\partial z_i}) dz_k
\wedge dz_i\,\, .$$ On the other hand, by definition of $f_i(z)$
($1 \le i \le n$) and definition of $\tilde{\theta}(z)$ from the
potential function $\Psi(z)$, we know that
$$f_i(z) = \sum_{\alpha =1}^{n-1}j_{\alpha}\tilde{\theta}_{\alpha}^{i}(z) +
m\tilde{\theta}_{n}^{i}$$
$$= \sum_{\alpha =1}^{n-1}j_{\alpha} \frac{\partial^2 \Psi}{\partial z_{\alpha} \partial z_k} + m \frac{\partial^2 \Psi}{\partial z_{n} \partial z_k}\,\, .$$
Since $\Psi(z)$ is holomorphic, it follows that
$$\frac{\partial}{\partial z_k}(\frac{\partial^2\Psi}{\partial z_{\alpha} \partial z_i}) =
\frac{\partial^3 \Psi}{\partial z_k \partial z_{\alpha} \partial z_i}
= \frac{\partial}{\partial z_i}(\frac{\partial^2\Psi}{\partial z_{\alpha} \partial z_k})\,\, .$$
Substituting this into the formula above, we obtain that $T_{\gamma}^* \omega_{\tilde{M}} = \omega_{\tilde{M}}$.
\end{proof}

\medskip
{\it (III) Group quotient of $\tilde{M}$ by $\Gamma = \Z^n$.}
\medskip

Let $M = \tilde{M}/\Gamma$. By Proposition \ref{discontinuous}, $M$
is a smooth symplectic manifold with symplectic form $\omega_{M}$ induced by $\omega_{\tilde{M}}$ and $M$ admits a fibration $f : M \to B$ induced by
$\tilde{f}$. We denote the (scheme theoretic) fiber $f^{-1}(b)$ over $b \in B$
by $M_b$. Let us describe the fibers $M_b$.

\medskip
{\it (III-1) Smooth fibers $M_b$}
\medskip

First consider the case where $z_n(b) \not= 0$, i.e., the case where $M_b$
is smooth. We have
$$\tilde{M}_b = (\C^{\times})^{n-1} \times \{(x_k, y_k) \vert x_ky_k = b_n\} \simeq
(\C^{\times})_{(w_1, \ldots, w_{n-1}, w_n)}^{n}\, \, .$$ and $w_n
= z_n(b)^k y_k$. Let $\langle e_i \rangle_{i=1}^n$ be the ordered
standard basis of $\Gamma$. From the description in (II), the
action of $\Gamma$ is given by:
$$T_{e_i}^* w_j = {\rm exp}(2\pi \sqrt{-1}\theta_i^j(b))w_j$$
$$T_{e_i}^* w_n =
{\rm exp}(2\pi \sqrt{-1}\theta_i^n(b))w_n$$ for $1 \le i \le n-1$
with $\theta_{i}^{j} = \tilde{\theta}_{i}^{j}$ and
$$T_{e_n}^* w_j = {\rm exp}(2\pi \sqrt{-1}\theta_n^j(b))w_j$$
$$T_{e_n}^* w_n =
{\rm exp}(2\pi \sqrt{-1}\tilde{\theta}_n^n(b))z_n(b)^{\ell}w_n =
{\rm exp}(2\pi \sqrt{-1}\theta_n^n(b))w_n\,\, .$$

Hence $M_b = \tilde{M}_b/\Gamma$ is an $n$-dimensional principally
polarized abelian variety of period $\theta(b)$, as desired. By
the description of $\omega_M$, the fibers $M_b,$ $z_n(b) \not= 0,$
are also Lagrangian submanifolds.

\begin{proposition}\label{relativepolarization} For $ b\in B, z_n(b) \neq 0$, choose the basis
$$p_{1,b}, \ldots , p_{n,b}\, q_{1,b}, \ldots , q_{n,b}$$
of $H_1(M_b, \Z)$ such that $\tilde{M}_b = \C^{n}/\langle p_{j,b} \rangle_{j=1}^{n}$ and
$$q_{i,b} = \sum_{j=1}^{n} \theta_i^j(b)p_{j,b}$$
for each $i$ ($1 \le i \le n$). Let
$$p_b^1, \ldots , p_b^n\, q_b^1, \ldots , q_b^n$$
be the dual basis of $H^1(M_b, \Z)$. Then the integral $2$-form
$$L_b := \sum_{i=1}^{n} p_b^i \wedge q_b^i$$
give a monodromy invariant principal polarization of $M$ over $B \setminus D$
where $D = (z_n= 0)$.
\end{proposition}
\begin{proof}
When $z_n(b) \neq 0$, the fiber $\tilde{M}_b$ of $\tilde{p}:
\tilde{M} \to B$ is $(\C^{\times})^n$ and this family has no
monodromy over $B \setminus (z_n =0)$. Thus we can fix a basis
$p_{1,b}, \ldots, p_{n,b}$ of $H_1(\tilde{M}_b, \Z)$ uniformly in
$b, z_n(b) \neq 0.$

To get a basis of $H_1(M_b, \Z),$ we  choose additional elements
$q_{1,b}, \ldots, q_{n,b} \in H_x(M_b, \Z)$ determined by the
deck-transformation of $\tilde{M}_b$ induced by the action
$T_{e_1}, \ldots, T_{e_n}$. From the description of $T_{e_i}^*$ on
$w_j$, they satisfy the relation $$q_{i,b} = \sum_{j=1}^{n}
\theta_i^j(b)p_{j,b}.$$ We see that $q_{1,b}, \ldots, q_{n-1, b}$
are invariant under the monodromy, while $q_{n,b} \mapsto q_{n,b}
+ \ell p_{n,b}$ under the monodromy of the generator $\gamma$ of
$\pi_1(B \setminus D)$, i.e., the circle around discriminant
divisor $z_n =0$. The 2-form $L_b$ is a principal polarization on
$M_b$. It remains to show that $L_b$ is invariant under the
monodromy.
  By definition of $\theta(b)$, we compute that
$$\gamma^{*}(L_b) =  \gamma^{*}(\sum_{i=1}^{n-1} p_{i,b} \wedge q_{i,b}) +
\gamma^{*}(p_{n, b} \wedge q_{n, b})$$
$$= \sum_{i=1}^{n-1} p_b^i \wedge q_b^i + p_b^n \wedge (q_b^n - \ell p_b^n) = L_b.$$
This implies the invariance.
\end{proof}

\begin{remark}\label{nakamura} As in \cite{Na}, one can also describe $\tilde{f} : \tilde{M} \to B$ in terms of toric geometry. Following  an argument similar to \cite{Na}, Section 4, it seems possible to give a relatively principally polarized divisor (the relative theta divsor) which is defined globally over $B \setminus D$.
However, its closure is not necessarily $f$-ample even if total
space is of dimension $4$ (cases of stable principally polarized
Lagrangian $4$-folds). In fact, a failure of $f$-ampleness of the
closure already happens when the fiber dimension $2$ and the base
dimension $1$ as explicitly described in \cite{Na} Section 4, Page
219. See also Remark 3.6.
\end{remark}

\medskip
{\it (III-2) Singular fibers $M_b$}
\medskip

Next consider the singular fibers of $f$. They are $M_{b}$ with $z_n(b) = 0$.
Recall from (I) that $\tilde{M}_b$ is the product of $(\C^{\times})^{n-1}$
with coordinate $(w_i)_{i=1}^{n-1}$ and the infinite tree $\cup_{k \in \Z} \BP_k^1$ of projective lines $\BP_k^1$ with affine coordinate $y_k$, and $M_b = \tilde{M}_b/\Gamma$.
Let us denote by $(0)_{k}, (\infty)_k \in \BP_k^1$ the two points on the projective line $\BP^1_k$ such that
$(0)_k$ is identified with $(\infty)_{k-1}$ in the tree.

From (II), the action of $\Gamma$ is given by:
$$T_{e_i}^* w_j = {\rm exp}(2\pi \sqrt{-1}\tilde{\theta}_i^j(b))w_j$$
$$T_{e_i}^* y_k =
{\rm exp}(2\pi \sqrt{-1}\tilde{\theta}_i^n(b))y_k$$
for $1 \le i \le n-1$ and
$$T_{e_n}^*w_j = {\rm exp}(2\pi \sqrt{-1}\tilde{\theta}_n^j(b))w_j$$
$$T_{e_n}^*y_k =
{\rm exp}(2\pi \sqrt{-1}\theta_n^n(b))y_{k-\ell}\,\, ,$$ where
$\langle e_i \rangle_{i=1}^n$ is the ordered standard basis of
$\Gamma$. Here we note that the last equality shows that the
monodromy operation corresponds to the shift of the components of
the infinite tree $\cup_{k \in \Z} \BP_k^1$. Thus
$\tilde{M}/<e_n>$ can be described as the variety obtained from
$$(\C^{\times})^{n-1} \times \cup_{k=0}^{\ell -1} \BP^1_k$$ by
identifying the point $$(w^1, \ldots, w^{n-1}) \times (0)_0 \; \in
\; (\C^{\times})^{n-1} \times (0)_0$$   with
 the point $$ (\exp ( 2\pi \sqrt{-1}  \theta^1_n) w^1, \ldots, \exp ( 2\pi \sqrt{-1}  \theta^{n-1}_n)w^{n-1})
  \; \in \; (\C^{\times})^{n-1} \times (\infty)_{\ell-1}.$$
 From this description, $M_b$ consists of $\ell$ irreducible
components, each of whose normalization is isomorphic to a
$\BP^1$-bundle over $(n-1)$-dimensional complex torus isogenous to
$(\C^{\times})^{n-1}/\langle e_i \rangle_{i=1}^{n-1}$, where the
action of $\langle e_i \rangle_{i=1}^{n-1}$ is given by the
coordinate action $T_{e_i}^*$ ($1 \le i \le n-1$) on $w_j$ ($1 \le
j \le n-1$) described above. Note that the quotient
$(\C^{\times})^{n-1}/\langle e_i \rangle_{i=1}^{n-1}$ is compact
because the imaginary part of $\tilde{\Theta}_1(z)$ is positive
definite from the assumption that the imaginary part of
$\tilde{\theta}$ is positive definite. We also note that the
characteristic cycles are of type $I_{m}$ for some $1 \le m \le
\infty$.

From this description,  the following is now clear:

\begin{theorem}\label{mainconst} The fibration $f : M \to B$ constructed above is a proper, flat,
 principally polarized stable Lagraingian fibration with a potential function $\Psi(z)$.
 Moreover, $\ell$ is the number of components of the singular fiber and $S^i_1 \neq S^i_2$ for all $i$ in the notation of Definition \ref{d.stable}.
\end{theorem}

Given a principally polarized Lagrangian fibration, we can find a
potential function $\Psi$ on $B$ as in Theorem \ref{t.main}.
Starting from $\Psi$ we can construct a principally polarized
stable Lagrangian fibration by Theorem \ref{mainconst}. These two Lagrangian fibrations
 must
agree outside the discriminant set. Thus they must be biholomorphic by the following.

\begin{proposition}\label{p.unique}
Let $f: M \to B$ and $f':M' \to B$ be two  Lagrangian fibrations
with the same discriminant $D \subset B$, having Lagrangian
sections $\Sigma \subset M$ and $\Sigma' \subset M'$.
Suppose there exists a biholomorphic morphism
$\Phi: M \setminus f^{-1}(D) \to M' \setminus f'^{-1}(D)$ such
that $\Phi(\Sigma) = \Sigma'$ and $\Phi$ is symplectomorphic, i.e.,
$\Phi^* \omega_{M'}= \omega_M$. Then $\Phi$ extends to a
biholomorphic morphism $M \cong M'$. \end{proposition}

\begin{proof} The proof is essentially given in the proof of Proposition 5.1 in \cite{HO2}.
 Let us sketch the argument. One can see that $\Phi$ is a bimeromorphic map between $M$ and $M'$.
Choose holomorphic coordinates $(z_1, \ldots, z_n)$ on $B$ such
that the discriminant $D$ is defined by $z_n=0$. Then the
Hamiltonian vector fields induced by $dz_1, \ldots, dz_{n-1}$
determine  $\C^{n-1}$-actions on $M$ and $M'$ such that $\Phi$ is
equivariant with respect to them. These $\C^{n-1}$-actions are
free and all orbits have dimension $n-1$. Since both $M$ and $M'$
have trivial canonical bundle, the exceptional loci of the
bimeromorphic map must be of codimension $\geq 2$. Since they are
invariant under the $\C^{n-1}$-actions, they must be union of
finitely many orbits of $\C^{n-1}$-actions. But then each
component of the exceptional loci in $M$ must be transformed to a
component of the exceptional loci in $M'$  biholomorphically. This
implies that there are no exceptional loci and $\Phi$ is
biholomorphic. \end{proof}

As a corollary of Theorem \ref{mainconst} and Proposition
\ref{p.unique}, we obtain

\begin{corollary}\label{c.number} For any
principally polarized stable Lagrangian fibration, the positive
integer $|\ell|$ in Theorem \ref{t.main} is the number of
components of the singular fiber and  $S^i_1 \neq S^i_2$ for each
$i$ in the notation of Definition \ref{d.stable}.
\end{corollary}

Theorem \ref{t.main}, Theorem \ref{mainconst} and Corollary \ref{c.number} complete the proof of
Theorem \ref{t.1}.

\section{Periods and the characteristic cycles}
\par
\noindent

In this section, we will examine the relation between types of the characteristic cycles and the periods.
For simplicity, we will restrict our discussion to the case of $\ell=1$. The generalization to arbitrary
  $\ell$ is straightforward. Explicit constructions (Constructions I -III) will be given when $n=2$, i.e.,
  constructions  of
$4$-dimensional principally polarized stable Lagrangian
fibrations.  Construction I gives an explicit example in which the
types of characteristic cycles change fiber by fiber. Construction
II gives an explicit example in which the types of characteristic
cycles are constant type $I_n$ ($n < \infty$) and Construction III
gives an explicit example in which the types of characteristic
cycles are constant type $A_{\infty}$.

\begin{proposition}\label{types} Let $f : M \to B$ be a $2n$-dimensional
principally polarized stable Lagrangian fibration with potential function
$\Psi(z)$ and $\ell =1$.
We denote the (univalent) period matrix by $\tilde{\theta}(z) =
(\tilde{\theta}_i^j(z))_{i, j=1}^{n}$ and the
multi-valued period matrix $\theta(z)$ of $f$ as
$$\theta(z) = \tilde{\theta}(z) +
\frac{\log z_n}{2\pi\sqrt{-1}}\left(\begin{array}{rr}
O_{n-1} & 0\\
0 & 1\end{array} \right)\,\, .$$ For $b  \in B$ for which $M_b$ is
singular, define $n(b)$ ($1 \le n(b) \le \infty$) to be the order
of
$$({\rm exp}(2\pi\sqrt{-1}\tilde{\theta}_n^j(b)))_{j=1}^{n-1}\, {\rm mod}\,
\langle ({\rm exp}(2\pi\sqrt{-1}\tilde{\theta}_i^j(b)))_{j=1}^{n-1}\, \vert\,
1 \le i \le n-1\, \rangle$$
in the multiplicative group
$$(\C^{\times})^{n-1}/\langle ({\rm exp}(2\pi\sqrt{-1}\tilde{\theta}_i^j(b)))_{j=1}^{n-1}\, \vert\,
1 \le i \le n-1\, \rangle\,\, .$$ Then the characteristic cycle of
$M_b$ is of type $I_{n(b)}$.
\end{proposition}

\begin{remark} The description above is simpler when $f : M \to B$ is a $4$-dimensional
principally polarized stable Lagrangian fibration with potential
function $\Psi(z) = \Psi(z_1, z_2)$ and $\ell =1$, as follows. We
shall use this description  in Constructions (I)-(III) below. We
write the (univalent) period matrix $\tilde{\theta}(z)$ and the
multi-valued period matrix $\theta(z)$ of $f$ as
$$\tilde{\theta}(z) = \left(\begin{array}{rr}
\tilde{\theta}_1(z) & \tilde{\theta}_2(z)\\
\tilde{\theta}_2(z) & \tilde{\theta}_3(z)
\end{array} \right)\,\, , \theta(z) = \tilde{\theta}(z) +
\frac{\log z_2}{2\pi\sqrt{-1}}\left(\begin{array}{rr}
0& 0\\
0 & 1\end{array} \right)\,\, .$$ For $b =(b_1, 0) \in B$ for which
$M_b$ is singular, the characteristic  cycle of $M_b$ is then of
type $I_{n(b)}$, where $n(b)$ ($1 \le n(b) \le \infty$) is exactly
the order of
$${\rm exp}(2\pi\sqrt{-1}\tilde{\theta}_2(b))\, {\rm mod}\,
\langle {\rm exp}(2\pi\sqrt{-1}\tilde{\theta}_1(b)) \rangle$$
in the multiplicative group $\C^{\times}/\langle {\rm exp}(2\pi\sqrt{-1}\tilde{\theta}_1(b)) \rangle$.
\end{remark}

\begin{proof} In the description (III-2), $M_b$ is the quotient of
$\tilde{M}_b = \cup_{k \in \Z} (\C^{\times})^{n-1} \times \BP_k^1$ with coordinates $((w_j)_{j=1}^{n-1}, y_k)$ ($k \in \Z$)
by the action of $\Gamma = \Z^n$ with ordered standard basis
$$\langle e_1, \ldots , e_{n-1} , e_{n}\rangle\,\, .$$ In terms of the standard basis, the action is given by:
$$T_{e_i}^* : (w_j)_{j=1}^{n-1} \mapsto
({\rm exp}(2\pi \sqrt{-1}\tilde{\theta}_i^j(b))w_j)_{j=1}^{n-1}$$
$$T_{e_i}^* : y_k \mapsto
{\rm exp}(2\pi \sqrt{-1}\tilde{\theta}_i^n(b))y_k$$
for $e=i$ ($1\le i \le n-1$), and for $e_n$
$$T_{e_n}^* : (w_j)_{j=1}^{n-1} \mapsto
({\rm exp}(2\pi \sqrt{-1}\tilde{\theta}_n^j(b))w_j)_{j=1}^{n-1}$$
$$T_{e_n}^* : y_k \mapsto
{\rm exp}(2\pi \sqrt{-1}\tilde{\theta}_n^n(b))y_{k-1}\,\, .$$ Thus
$\tilde{M}_b/\langle e_n \rangle$ is $(\C^{\times})^{n-1} \times
\BP_0^1$ in which $((w_j)_{j=1}^{n-1}, 0)$ and
$((w_j')_{j=1}^{n-1}, \infty)$ are  identified exactly when the
two points $(w_j)_{j=1}^{n-1}$ and $(w_j')_{j=1}^{n-1}$ of
$(\C^{\times})^{n-1}$ are in the same orbit under the action of
the cyclic subgroup
$$G(b) := \langle ({\rm exp}(2\pi \sqrt{-1}\tilde{\theta}_n^j(b)))_{j=1}^{n-1} \rangle$$
of $({\mathbf C}^{\times})^{n-1}$. On the other hand,
$((\C^{\times})^{n-1} \times \BP_0^1)/\langle e_i
\rangle_{i=1}^{n-1}$ is the normalization of $M_b$. Thus, $M_b$ is
obtained from $((\C^{\times})^{n-1} \times \BP_0^1)/\langle e_i
\rangle_{i=1}^{n-1}$ by identifying the two $(n-1)-$dimensional
complex tori $((\C^{\times})^{n-1} \times \{\infty\})/\langle e_i
\rangle_{i=1}^{n-1}$ and $((\C^{\times})^{n-1} \times \{0
\})/\langle e_i \rangle_{i=1}^{n-1}$,by the action of $G(b)$
above. Here, as a subgroup of $(\C^{\times})^{n-1}$, the group
$\langle e_i \rangle_{i=1}^{n-1}$ is the multiplicative subgroup
generated by the $n-1$ elements
$$({\rm exp}(2\pi \sqrt{-1} \tilde{\theta}_i^j(b)))_{j=1}^{n-1}\,\, ,\,\,
 1 \le i \le n-1\,\, .$$ This implies the result.
\end{proof}

\medskip
{\it Construction I.}
\medskip

Under the notation of Section 5, we set $\ell = 1$ and
$$\Psi(z_1, z_2) := \frac{(z_1 + 5\sqrt{-1})^3 + (z_2 + 5\sqrt{-1})^3 +
3z_1^2z_2 + 3z_1z_2^2}{6}\, ,$$
Then
$$\tilde{\theta}(z) = \left(\begin{array}{rr}
z_1 + z_2 + 5\sqrt{-1} & z_1+z_2\\
z_1+z_2 & z_1 + z_2 + 5\sqrt{-1}
\end{array} \right)\,\, , \theta(z) = \tilde{\theta}(z) +
\frac{\log z_2}{2\pi\sqrt{-1}}\left(\begin{array}{rr}
0& 0\\
0 & 1\end{array} \right)\,\,$$ and
$${\rm Im}\, \tilde{\theta}(z) =
\left(\begin{array}{rr}
y_1 + y_2 + 5 & y_1+y_2\\
y_1 +y_2 & y_1 + y_2 + 5
\end{array} \right)\,\, .$$
Here and hereafter $x_i$ and $y_i$ are the real and imaginary
part of $z_i$ respectively. Since $t+5 > 0$ and $(t+5)^2 - t^2 >
0$ when $-2 < t < 2$, it follows that ${\rm Im}\, \theta(z)$ is
positive definite on the polydisk
$$\{(z_1, z_2)\, \vert\, \vert z_i \vert < 1\, \}\, .$$
Taking a smaller $2$-dimensional polydisk $B$ with multi-radius $\epsilon$,
we then obtain a $4$-dimensional Lagrangian
fibration $f : M \to B$,
associated with the potential function $\Psi(z)$ and $\ell = 1$.
The discriminant set is $z_2 = 0$. Define $N = N(z)$ to be the order of
$$e^{2\pi \sqrt{-1}z_1}\, {\rm mod}\, \langle
e^{2\pi \sqrt{-1}(z_1+ 5\sqrt{-1})} \rangle$$
in the multiplicative group $\C^{\times}/\langle
e^{2\pi \sqrt{-1}(z_1+ 5\sqrt{-1})} \rangle$. By abuse of language, we
include $N = \infty$ when the order is not finite. Then,
the characteristic cycle
of $M_{(z_1, 0)}$ is of type $I_N$.

\begin{proposition}\label{moving infty} In Construction 1, the characteristic cycle on $M_{(z_1, 0)}$
is of Type $I_k$ with $k < \infty$ if and only if
$$z_1 \in \Q(\sqrt{-1}) \,\, .$$
So, the singular fibers of finite characteristic cycle $I_k$ ($k <
\infty$) and the singular fibers of infinite characteristic cycle
$I_{\infty}$ are both dense over the disciminant set.  Moreover,
the characteristic cycle of $M_{(z_1, 0)}$ is precisely of type
$I_k$ ($k < \infty$) for $z_1 = 1/k$. So, the singular fibers with
characteristic cycles of type $I_k$ with any sufficiently large
$k$ appear in this family.
\end{proposition}

\begin{proof} By the definition of $N = N(z)$,
it follows that $N < \infty$ for $M_{(z_1, 0)}$ if and only if
there are integers $k > 0$ and $m$ such that
$$(e^{2\pi \sqrt{-1}z_1})^k = (e^{2\pi \sqrt{-1}(z_1+5\sqrt{-1})})^k\,\, .$$
The last condition is equivalent to
$$kz_1 - m(z_1 + 5\sqrt{-1}) \in \Z$$
which is also equivalent to
$$(k-m)x_1 \in \Z\,\, {\rm and}\,\, (k - m)y_1 - 5 m = 0\,\, .$$
Note that $k -m \not= 0$ in the last equivalent condition, as
otherwise $k = m = 0$. It is immediate to see that two integers $k
> 0$ and $m$ satisfying  last equivalent condition exist if and
only if $z_1 \in \Q(\sqrt{-1})$. Since $\vert e^{2\pi
\sqrt{-1}(z_1+5\sqrt{-1})} \vert > 1$ for $\vert z_1 \vert < 1$,
whereas $\vert e^{2\pi \sqrt{-1}/k} \vert = 1$ for $k \in \Z$, it
follows that the order $N(z)$ for $z_1 = 1/k$ is precisely the
order of $e^{2\pi \sqrt{-1}/k}$ in the multiplicative group
$\C^{\times}$. This implies the last statement.
\end{proof}

\medskip
{\it Construction II.}
\medskip

Under the notation of Section 5, we set $\ell = 1$ and
$$\Psi(z_1, z_2) := \frac{\sqrt{-1}(z_1^2 + z_2^2)}{2} +
\frac{z_1z_2}{k}\, ,$$
where $n$ is a positive integer.
Then
$$\tilde{\theta}(z) = \left(\begin{array}{rr}
 \sqrt{-1} & 1/k\\
1/k & \sqrt{-1}
\end{array} \right)\,\, , \theta(z) = \tilde{\theta}(z) +
\frac{\log z_2}{2\pi\sqrt{-1}}\left(\begin{array}{rr}
0& 0\\
0 & 1\end{array} \right)\,\,$$ and
$${\rm Im}\, \tilde{\theta}(z) =
\left(\begin{array}{rr}
1 & 0\\
0 & 1
\end{array} \right)\,\, .$$
The matrix ${\rm Im}\, \tilde{\theta}(z)$ is positive definite.
So, taking a smaller $2$-dimensional polydisk $B$,
we obtain $4$-dimensional Lagrangian
fibration $f : M \to B$,
associated with the potential function $\Psi(z)$ and $\ell = 1$ above.
The discriminant set is $z_2 = 0$. The order of
$$e^{2\pi \sqrt{-1}/k}\, {\rm mod}\, \langle
e^{-2\pi} \rangle$$
is exactly $k$
in the multiplicative group $\C^{\times}/\langle
e^{-2\pi} \rangle$,
where $-2\pi = 2\pi \sqrt{-1} \cdot \sqrt{-1}$. Then,
the characteristic cycle of $M_{(z_1, 0)}$ is of type $I_k$, and in particular, the type is constant.

\medskip
{\it Construction III.}
\medskip

Under the notation of Section 5, we set $\ell = 1$ and
$$\Psi(z_1, z_2) := \frac{\sqrt{-1}(z_1^2 + z_2^2)}{2} +
\alpha z_1z_2\, ,$$
where $\alpha$ is any irrational, real number, say $\sqrt{2}$.
Then
$$\tilde{\theta}(z) = \left(\begin{array}{rr}
 \sqrt{-1} & \alpha\\
\alpha & \sqrt{-1}
\end{array} \right)\,\, , \theta(z) = \tilde{\theta}(z) +
\frac{\log z_2}{2\pi\sqrt{-1}}\left(\begin{array}{rr}
0& 0\\
0 & 1\end{array} \right)\,\,$$ and
$${\rm Im}\, \tilde{\theta}(z) =
\left(\begin{array}{rr}
1 & 0\\
0 & 1
\end{array} \right)\,\, .$$
The matrix ${\rm Im}\, \theta(z)$ is positive definite.
So, taking a smaller $2$-dimensional polydisk $B$,
we obtain $4$-dimensional Lagrangian
fibration $f : M \to B$,
associated with the potential function $\Psi(z)$ and $\ell = 1$ above.
The discriminant set is $z_2 = 0$. Since $\alpha$ is irrational real number,
the element
$$e^{2\pi \sqrt{-1}\cdot \alpha}\, {\rm mod}\, \langle
e^{-2\pi} \rangle$$
is of infinite order
in the multiplicative group $\C^{\times}/\langle
e^{-2\pi} \rangle$,
where $-2\pi = 2\pi \sqrt{-1} \cdot \sqrt{-1}$. Then,
the characteristic cycle of $M_{(z_1, 0)}$ is of type $A_{\infty}$,
and in particular the type is constant.

\end{document}